\documentclass[twoside]{article}

\usepackage[margin=1in]{geometry}

\title{GPU-Accelerated Particle Methods for Evaluation of Sparse Observations for PDE-Constrained Inverse Problems}%
\author{Jeff Borggaard, Nathan E. Glatt-Holtz, Justin A. Krometis\\
  \scriptsize{emails: jborggaard@vt.edu, negh@tulane.edu, jkrometis@vt.edu}}
\date{}

\usepackage{framed,multirow}

\usepackage{amssymb}
\usepackage{latexsym}

\usepackage{url}
\usepackage{xcolor}
\definecolor{newcolor}{rgb}{.8,.349,.1}

\usepackage{amsfonts, amssymb, amsmath, amsthm, multicol,bbm}
\usepackage[colorlinks=true, pdfstartview=FitV, linkcolor=blue,
            citecolor=blue, urlcolor=blue]{hyperref}
\usepackage{graphicx}
\usepackage[colorinlistoftodos]{todonotes}

\usepackage[capitalize,nameinlink,noabbrev]{cleveref} 
\usepackage[section]{algorithm}      
\usepackage{algpseudocode}    

\usepackage{mathtools} 
\usepackage{xparse}
\usepackage{mathrsfs}

\newcommand{\R}{\mathbb{R}}           
\newcommand{\Z}{\mathbb{Z}}           

\DeclareMathOperator*{\argmin}{arg\,min}

\newcommand{\Ito}{It\^{o}~}

\newcommand{\Frechet}{Fr\'echet~}

\DeclarePairedDelimiterX\innerp[2]{\langle}{\rangle}{#1,#2}
\newcommand{\ip}[2]{ \innerp{#1}{#2} }

\newcommand{\norm}[1]{\left\lVert#1\right\rVert}

\NewDocumentCommand\Exp{g}{%
    \mathbb{E}\IfNoValueTF{#1}{}{\left[ #1 \right]}%
}

\NewDocumentCommand\Prob{g}{%
    \mathbb{P}\IfNoValueTF{#1}{}{\left[ #1 \right]}%
}

\newcommand{\unknown}{u}
\newcommand{\x}{\mathbf{x}}
\newcommand{\e}{\mathbf{e}}
\newcommand{\G}{\mathcal{G}}
\newcommand{\Obs}{\mathcal{O}}
\newcommand{\Sol}{\mathcal{S}}
\newcommand{\kbf}{\mathbf{k}}

\newcommand{\spatdom}{D}
\newcommand{\conductivity}{\kappa}
\newcommand{\pdesol}{\theta}
\newcommand{\pdesolvec}{\Theta}
\newcommand{\vfieldnobf}{v}
\newcommand{\vfield}{\mathbf{\vfieldnobf}}

\newcommand{\data}{\mathcal{Y}}

\newcommand{\brownm}{\mathbf{W}}
\newcommand{\X}{\mathbf{X}}

\newcommand{\bigO}{O} 

\newcommand{\bcf}{\pdesol_{\partial\spatdom}}
\newcommand{\stoptime}{\tau_{\spatdom}}

\newtheorem{Thm}{Theorem}[section]

\newtheorem{Prb}[Thm]{Problem}

\begin{document}

\markboth{J. Borggaard, N. Glatt-Holtz, J. Krometis}
{GPU-Accelerated Particle Methods for Evaluation of Sparse Observations for PDE-Constrained Inverse Problems}%

\maketitle

\begin{abstract}
We consider the inverse problem of estimating parameters of a driven diffusion (e.g., the underlying fluid flow, diffusion coefficient, or source terms) from point measurements of a passive scalar (e.g., the concentration of a pollutant). We present two particle methods that leverage the structure of the inverse problem to enable efficient computation of the forward map, one for time evolution problems and one for a Dirichlet boundary-value problem. The methods scale in a natural fashion to modern computational architectures, enabling substantial speedup for applications involving sparse observations and high-dimensional unknowns. Numerical examples of applications to Bayesian inference and numerical optimization are provided.
\end{abstract}

{\noindent \small {\it \bf Keywords:} Inverse Problems, Optimization, Scientific Computing, Parallel Computing, Passive Scalars}

\setcounter{tocdepth}{1}
\tableofcontents

\newpage

\section{Introduction}

Much recent computational and theoretical work has been devoted to the inverse problem of the estimation of unknown or optimal model parameters from finite observations of the model output \cite{kaipio2005statistical,dashti2017bayesian,hinze2009optimization,biegler2003large}.
For example, recent works
included reconstruction of a seismic wave speed field from waves recorded at a finite number of point receivers \cite{bui2013computational}, estimation of an ice sliding coefficient field from finite velocity observations \cite{petra2014computational}, and determination of the source of a chemical/biological attack from measurements of toxins \cite{boggs2006rapid}.

Approaches to these inverse problems typically involve many evaluations of the model, also called the forward map. For example, each step of a Markov Chain Monte Carlo (MCMC) method \cite{brooks2011handbook} for Bayesian inference or iteration of a numerical optimization routine \cite{nocedal2006numerical} will require evaluation of the forward map. When that map includes observations of the solution of a partial differential equation (PDE), these methods in turn require many solutions of the PDE. In many applications, these PDE solves are computationally-intensive and dominate the overall time to solution. 
For applications where the observations are sparse, i.e., low-dimensional, solving the PDE can be wasteful in the sense that a high-dimensional quantity is computed only to take a low-dimensional projection of it. It is therefore desirable to identify methods that allow computation of the parameter-to-observation map directly, without a full PDE solve. 

In this work, we consider the problem of estimating the parameters of a driven diffusion (e.g., the background flow, diffusion coefficient, or sources) from point measurements of a passive scalar (e.g., the concentration of a contaminant). We present two numerical methods that can be used to compute observations for PDE-constrained inverse problems without computation of the full scalar field. These methods therefore bypass the need to approximate a high-dimensional PDE solution at each step of the inverse problem and instead replace the full PDE solve with an array of particle solutions that are much less computationally expensive. Moreover, since the particle simulations are decoupled, they can be parallelized in a straightforward manner on modern computational architectures. The result is a dramatic speedup, particularly for problems in which the dimension of the unknowns is significantly larger than the dimension of the observations.

\section{Motivation}

In this article, we consider two inverse problems: 

\begin{Prb}[Time-dependent Advection-Diffusion]
  Let $\spatdom$ be an open, bounded subset of $\R^n$. Let $\vfield:\spatdom \to \R^n$, $\sigma:\spatdom \to \R^{n \times n}$, and $\pdesol_0:\spatdom \to \R$ be given functions. Assume that there exists a $\pdesol:\overline{\spatdom} \to \R$ satisfying
  \begin{equation}
    \frac{d\pdesol}{dt}(\x) = -\vfield(\x) \cdot \nabla \pdesol(t,\x) + \frac{1}{2} \sum_{i,j} (\sigma\sigma^T)_{i,j}(\x) \frac{\partial^2}{\partial x_i \partial x_j} \pdesol(t,\x), \quad \pdesol(0,\x)=\pdesol_0(\x).
    \label{eq:pde_diffusion}
  \end{equation}
  \textbf{Goal:} Estimate unknown $\unknown$ (e.g., one or more of: advection field $\vfield$, diffusion coefficients $\sigma$, or initial condition $\pdesol_0$) from finite, possibly noisy point observations $\pdesol(t_j,\x_j;\unknown)$, $j=1,\dots,N$.
  \label{prob:ad}
\end{Prb}

\begin{Prb}[Boundary Value (Dirichlet) Problem]
  Let $\spatdom$ be an open, bounded subset of $\R^n$. Let $\vfield:\spatdom \to \R^n$, $\sigma:\spatdom \to \R^{n \times n}$, $f:\spatdom \to \R$, and $\bcf:\partial \spatdom \to \R$ be given functions. Assume that there exists a $\pdesol:\overline{\spatdom} \to \R$ satisfying
  \begin{equation}
    \begin{aligned}
      -\vfield(\x) \cdot \nabla \pdesol(\x) + \frac{1}{2} \sum_{i,j} (\sigma\sigma^T)_{i,j}(\x) \frac{\partial^2}{\partial x_i \partial x_j} \pdesol(\x) &= f(\x) && \x \in \spatdom \\
      \pdesol(\x) &= \bcf(\x) && \x\in \partial\spatdom.
    \end{aligned}
    \label{eq:bvp}
  \end{equation}
  \textbf{Goal:} Estimate unknown $\unknown$ (e.g., one or more of: advection field $\vfield$, diffusion coefficients $\sigma$, forcing $f$, or boundary condition $\bcf$) from finite point observations $\pdesol(\x_j;\unknown)$, $j=1,\dots,N$.
  \label{prob:bvp}
\end{Prb}

In each problem, the goal is to estimate the parameters of a driven diffusion from point measurements of the passive scalar $\pdesol$. Solving these inverse problems in practice will typically involve many computations of the forward map from a given parameter $\unknown$ to its associated observations, which we denote by 
\begin{equation}
  \G(\unknown)=\left\{ \G_j(\unknown) \right\}_{j=1}^N,\text{ where }
  \begin{cases} 
    \G_j(\unknown) = \pdesol(t_j,\x_j,\unknown) &\text{ for \cref{prob:ad}} \\
    \G_j(\unknown) = \pdesol(\x_j,\unknown)     &\text{ for \cref{prob:bvp}}.
  \end{cases}
  \label{eq:Gmap}
\end{equation}
$\G$ can be thought of as being composed of two operators
\begin{equation}
	\G(\unknown) = \Obs \circ \Sol(\unknown)
\end{equation}
where
\begin{itemize}
  \item The solution operator $\Sol:\unknown \mapsto \pdesol$ requires solving the PDE \eqref{eq:pde_diffusion} or \eqref{eq:bvp} for parameter $\unknown$
  \item The observation operator $\Obs$ involves taking point observations from $\pdesol(\unknown)$
\end{itemize}
As such, $\G$ is typically evaluated in two steps:
\begin{enumerate}
  \item Compute $\pdesol(\unknown)$ via some numerical PDE solver
  \item Compute observations $\G(\unknown)$ from $\pdesol(\unknown)$
\end{enumerate}

This natural approach has the benefit of allowing the application of third-party ``black-box'' PDE solvers to the inverse problem. However, computing $\Sol(\unknown)$ involves approximating a solution that is infinite-dimensional, which can be very computationally expensive, requiring a discretization with many thousands or millions of degrees of freedom. By contrast, the evaluation of $\Obs$ then involves projecting that PDE solution into a finite-dimensional space. As a result, much of the work involved in approximating $\Sol$ is, in some sense, discarded in the application of $\Obs$. In the next section we present a numerical method for evaluation of $\G$ that breaks this two-step paradigm.

\section{Methods}
In this section, we present a particle method that will allow point evaluation of $\pdesol$ directly from the unknown $\unknown$ without separate approximation of $\Sol$. The method will involve simulating an ensemble of particles (\Ito diffusions), which is a well-known method for approximating $\pdesol$; see, for example, \cite{chorin1973numerical}, \cite{degond1989weighted}, or \cite{majda2002vorticity} for details. However, to leverage the sparse nature of the observations, for this application we will simulate the particles \emph{backward in time} from their final condition to their initial condition. Doing so will allow us to avoid computing the entire field $\pdesol$ by computing it only where it is needed. 

\subsection{A Particle Method for \cref{prob:ad}}\label{sec:ad_method}

In this section, we present a particle method for computing $\G(\unknown)$ for the time-dependent problem \cref{prob:ad}. A key ingredient is Kolmogorov's Backward Equation, which is presented in \cref{thm:kolmogorov_backward}. See, for example, Section 8.1 of \cite{oksendal2013stochastic} for details.

\begin{Thm}[Kolmogorov's Backward Equation]
  Suppose that for all $t>0$ and $\x \in \spatdom$, $\pdesol$ satisfies \eqref{eq:pde_diffusion} 
  with $\pdesol_0 \in C^2(\spatdom)$ and $\vfield,\sigma \in C^1(\spatdom)$. Then $\pdesol(t,\x)$ is given by
  \begin{equation}
    \pdesol(t,\x)=\Exp \pdesol_0(\X_t)
    \label{eq:kolmogorov_backward}
  \end{equation}
  where $\X_t$ is the \Ito diffusion
  \begin{equation}
    d\X_t = -\vfield(\X_t)dt + \sigma(\X_t) d\brownm_t, \quad \X_0=\x.
    \label{eq:ito_diffusion}
  \end{equation}
  \label{thm:kolmogorov_backward}
\end{Thm}

\begin{figure}[htbp]
  \centering
  \includegraphics[width=0.4\textwidth]{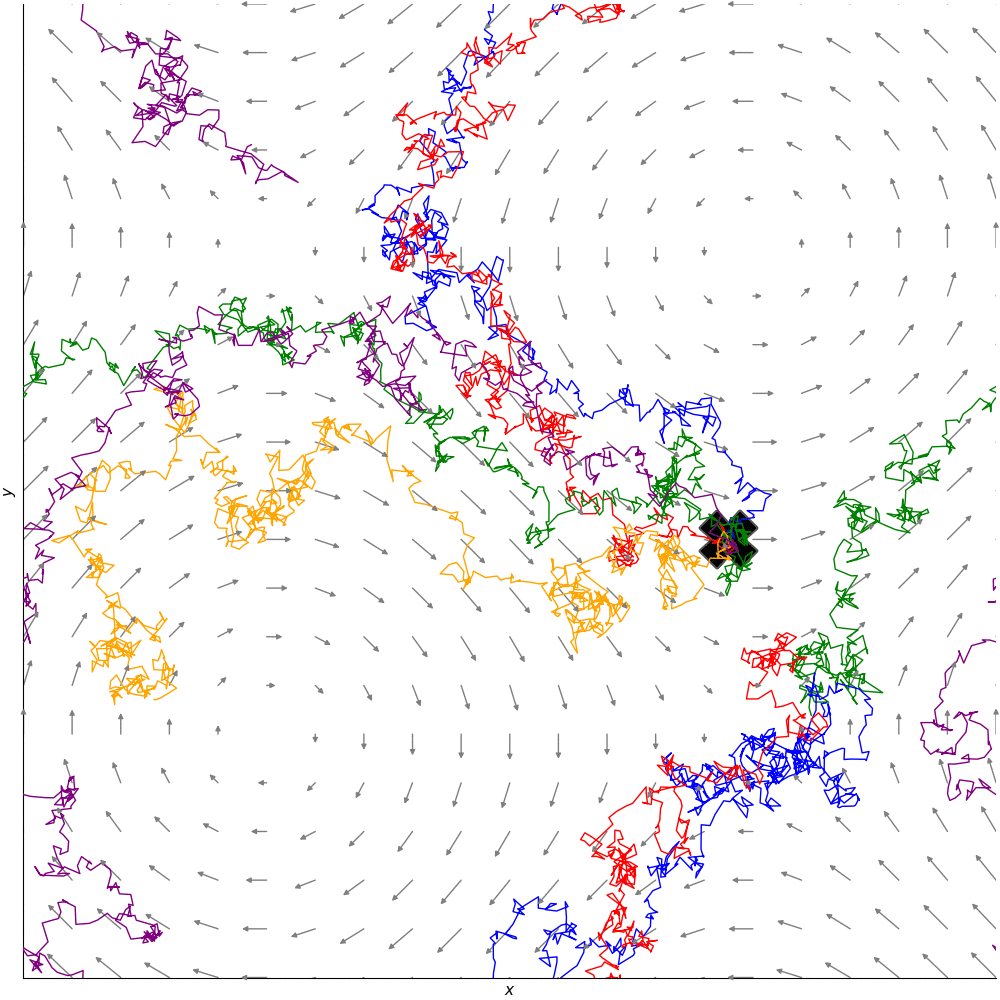}
  \caption{Traces of 5 simulated \Ito diffusions with the same final position.}
  \label{fig:particle_trace}
\end{figure}

Kolmogorov's Backward Equation tells us that the value of $\pdesol$ at a particular time and location $(t,\x)$ is given by the average value of the initial condition evaluated at the position of the \Ito diffusion \eqref{eq:ito_diffusion} at time $t$ when initialized at $\x$. This suggests a numerical method for evaluating $\G_j(\unknown)=\pdesol(t_j,\x_j;\unknown)$: (1) initialize a series of particles from $\x_j$; (2) simulate their movement to time $t_j$ according to \eqref{eq:ito_diffusion}; (3) evaluate $\pdesol_0$ at that location; and (4) take the average. See \cref{fig:particle_trace} for an illustration.

Numerical integration of \eqref{eq:ito_diffusion} can be computed, for example, with an Euler-Maruyama approximation \cite{kloeden1999numerical,higham2001ain}
\begin{equation}
  \X_{i+1} = \X_i - \vfield(\X_i)\Delta t_{i} + \sigma(\X_i)\sqrt{\Delta t_{i}}\xi_i,
  \label{eq:euler}
\end{equation}
or a Milstein approximation
\begin{equation}
  \X_{i+1} = \X_i - \vfield(\X_i)\Delta t_{i} + \sigma(\X_i)\sqrt{\Delta t_{i}}\xi_i + \frac{1}{2}\sigma(\X_i)\sigma'(\X_i)\left( \xi_i^2 - 1 \right)\Delta t_{i},
  \label{eq:millstein}
\end{equation}
where $\X_{i} = \X(t_i)$, $\Delta t_{i} = t_{i+1}-t_{i}$, and $\xi_i \sim N(0,1)$. The resulting algorithm is described in \cref{alg:ad_method}, where $N_{o}$ is the number of observations and $N_p$ is the number of particles used per observation. 

\begin{algorithm}
\caption{Particle Method for Computing $\G_j(\unknown)$.}\label{alg:ad_method}
\begin{algorithmic}[1]
  \State Given $\unknown$ (e.g., $\vfield$, $\sigma$, and/or $\pdesol_0$)
  \For{$j=1 \dots N_o$}
    \For{$i=1 \dots N_p$}
      \State Set $\X_0^{(i)} = \x_j$
      \State Simulate \eqref{eq:ito_diffusion} with \eqref{eq:euler} or \eqref{eq:millstein} to get $\X_{t_j}^{(i)}(\unknown)$
      \State Compute $g_j^{(i)}(\unknown)=\pdesol_0\left(\X_{t_j}^{(i)}\right)$
    \EndFor
    \State Compute $\G_j(\unknown)=\pdesol(t_j,\x_j) \approx \frac{1}{N_p}\sum_{i} g_j^{(i)}$
  \EndFor
\end{algorithmic}
\end{algorithm}

\subsection{A Particle Method for \cref{prob:bvp}}\label{sec:bvp_method}
In this section, we present a particle method for computing $\G(\unknown)$ for the boundary value problem \cref{prob:bvp}. The key idea is in the following theorem, which identifies the solution to \eqref{eq:bvp} in terms of \Ito diffusions; see, e.g., Section 9.1 of \cite{oksendal2013stochastic} for details.
\begin{Thm}[Particle Solution to \eqref{eq:bvp}]
  Let $f \in C(\spatdom)$, $\bcf \in C(\partial \spatdom)$, $\vfield \in C^1(\spatdom)$, and $\sigma \in C^1(\spatdom)$. Suppose there exists $\pdesol \in C^2(\spatdom)$ satisfying \eqref{eq:bvp}. Then $\pdesol$ is given by
  \begin{equation}
    \pdesol(\x) = \Exp \left[ \bcf\left(\X_{\stoptime}\right) \chi_{\stoptime < \infty} \right] - \Exp \left[ \int_0^{\stoptime} f(\X_t) \,dt \right]
    \label{eq:bvp_solution}
  \end{equation}
  where, as in \eqref{eq:ito_diffusion}, $\X_t$ is the \Ito diffusion
  \begin{equation*}
    d\X_t = -\vfield(\X_t)dt + \sigma(\X_t) d\brownm_t, \quad \X_0=\x
  \end{equation*}
  and $\stoptime$ denotes the first exit time from $\spatdom$.
\end{Thm}

As in \cref{sec:ad_method}, this theorem presents a natural numerical approach to evaluating $\G_j(\unknown)=\pdesol(\x_j,\unknown)$: (1) initialize a series of particles from $\x_j$; (2) simulate their movement according to \eqref{eq:ito_diffusion} until they exit $\spatdom$, integrating $f(\X_t) \,dt$ along the way; (3) evaluate $\bcf$ at the boundary location; (4) subtract the time integral of $f(\X_t)$; and (5) take the average. The resulting algorithm is described in \cref{alg:bvp_method}, where $N_{o}$ is the number of observations and $N_p$ is the number of particles used per observation.
\begin{algorithm}
  \caption{Particle Method for Computing $\G_j(\unknown)$ for \cref{prob:bvp}.}\label{alg:bvp_method}
\begin{algorithmic}[1]
\State Given $\unknown$ (e.g., $\vfield$, $\sigma$, $f$, $\bcf$)
  \For{$j=1 \dots N_o$}
    \For{$i=1 \dots N_p$}
      \State Set $\X_0^{(i)} = \x_j$
      \While{$\X_t^{(i)} \in \spatdom$}
        \State Simulate particle position $\X_t^{(i)}$ with \eqref{eq:euler} or \eqref{eq:millstein}
      \EndWhile
      \State Compute $g_j^{(i)}(\unknown)=\bcf\left(\X_{\stoptime^{(i)}}^{(i)}\right) - \int_0^{\stoptime^{(i)}} f(\X_t^{(i)}) \,dt$
    \EndFor
    \State Compute $\G_j(\unknown)=\pdesol(t_j,\x_j) \approx \frac{1}{N_p}\sum_{i} g_j^{(i)}$
  \EndFor
\end{algorithmic}
\end{algorithm}

\subsection{Parallelization}\label{sec:parallel}
Note that for a given observation point, all but the final step of \cref{alg:ad_method} or \cref{alg:bvp_method} can be computed in parallel. In addition, the steps for separate observation points are entirely independent. As a result, the algorithms are embarrassingly parallel and can therefore be parallelized in a straightforward manner using any number of computational paradigms, such as message passing interface (MPI) processes or OpenMP threads, and naturally vectorize to leverage ``single instruction multiple data'' (SIMD) capabilities on modern CPUs or GPUs. 

\begin{figure}[htbp]
  \centering
  \includegraphics[width=0.3\textwidth]{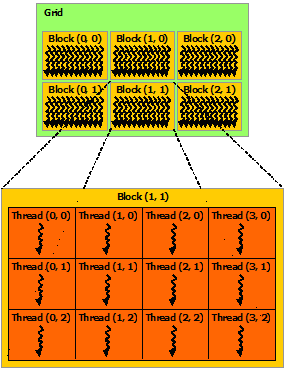}
  \caption{GPU thread blocks \cite{nvidia2018nvidia}.}
  \label{fig:cuda_thread_blocks}
\end{figure}

For example, \cref{fig:cuda_thread_blocks}, from NVIDIA's documentation \cite{nvidia2018nvidia}, 
illustrates the layout of threads and blocks on an NVIDIA GPU. Each thread is a single execution unit that is grouped into a block and then run across one or more SIMD. \cref{alg:ad_method} or \cref{alg:bvp_method} can be ported to this architecture in a natural fashion: each observation can be assigned to a block or group of blocks, with each particle run in a separate thread. A single GPU can execute thousands of threads at a time, allowing thousands of particles to be simulated simultaneously on a single chip. Moreover, larger problems can be spread across multiple GPUs or multiple machines using MPI to gain even greater efficiency; the overhead is minimal as only the average value for each observation needs to be returned to the master process.

\subsection{Computational Complexity}
In this section, we will compare the computational costs of \cref{alg:ad_method} to that of a reference Galerkin-based PDE solve. The analysis for \cref{alg:bvp_method} is similar and we summarize the results at the end of this section. Consider the case where $\unknown$ is approximated by a basis expansion $\unknown \approx \sum_{i=1}^{N_\unknown} \vfieldnobf_i \e_i$ and evaluation of each basis function $\e_i$ has computational cost $C_{b}$. We use the Euler-Maruyama approximation \eqref{eq:euler} and assume $N_t$ timesteps per observation. Further, we assume evaluating $\G$ requires $N_{o}$ observations and use $N_p$ particles per observation. Then the computational cost of \cref{alg:ad_method} for $P$ parallel processes/threads is given by 
\begin{equation}
  C_{particle} = \bigO \left( \frac{1}{P} N_{o} N_{p} N_{t} N_{\unknown} C_{b} \right).
  \label{eq:cost_particle}
\end{equation}
Since the observations and particle simulations are almost entirely independent, they can be executed in parallel (see \cref{sec:parallel}). So for large $P$ we have the limit
\begin{equation}
  C_{particle} \to \bigO \left( N_{t} N_{\unknown} C_{b} \right).
  \label{eq:cost_particle_limit}
\end{equation}

A more traditional method of solving the PDE would be to use a Galerkin projection, which would involve projecting the PDE \eqref{eq:pde_diffusion} onto a set of basis functions $\{ \phi \}_{l=1}^{N_b}$ to get a system of ODEs for the coefficients of $\pdesol$:
\begin{equation}
  \begin{aligned}
    M\dot{\pdesolvec} &= A\pdesolvec \\
    M_{lm} &= \ip{\phi_l}{\phi_m} \\
    A_{lm} &= \left\langle \phi_l,-\vfield(\x) \cdot \nabla \phi_m + \frac{1}{2} \sum_{i,j} (\sigma\sigma^T)_{i,j}(\x) \frac{\partial^2}{\partial x_i \partial x_j} \phi_m \right\rangle. \\
  \end{aligned}
  \label{eq:galerkin}
\end{equation}
The bases $\phi_{l}$ could, for example, be Fourier or finite element basis functions. The system \eqref{eq:galerkin} is then integrated by repeated iteration of some combination of $M,A$ (explicit, Runge-Kutta methods) and/or their inverses (implicit methods). \cref{alg:reference} outlines this algorithm for Explicit Euler time integration.
\begin{algorithm}
\caption{Reference Method for Computing $\G_j(\unknown)$ (Galerkin, Explicit Euler).}\label{alg:reference}
\begin{algorithmic}[1]
  \State Given $\unknown$ (e.g., $\vfield$, $\sigma$, and/or $\pdesol_0$)
  \State Project \eqref{eq:pde_diffusion} onto $\{ \phi \}_{l=1}^{N_b}$ to get system $M\dot{\pdesolvec} = A\pdesolvec$
  \For{$i=1 \dots N_t$}
    \State Multiply $\pdesolvec_{i} = \left( M + \Delta t_i A \right)\pdesolvec_{i-1}$
  \EndFor
  \For{$j=1 \dots N_o$}
    \State Compute $\G_j(\unknown)=\sum_{l} \pdesolvec_{l}(t_j)\phi_l(\x_j)$
  \EndFor
\end{algorithmic}
\end{algorithm}

The computational costs of assembling the matrix $A$ are heavily dependent on the choice of $\vfield$ and $\{\phi\}$ and therefore difficult to characterize in general. The cost of computing the observations is typically small. Ignoring these two factors, the computational cost of \cref{alg:reference} is dominated by the time integration of the system, which is made up of a series of matrix-vector multiplications. In general, the cost of this computation for $N_b$ basis functions and $N_t$ time steps is
\begin{equation}
  C_{reference} \to \bigO \left( N_{t} N_{b}^2 \right).
  \label{eq:cost_reference}
\end{equation}
Note also that to model an unknown of dimension $N_{\unknown}$ requires $N_b \ge N_{\unknown}$. In practice, however, accurate modeling of $\pdesol$ may require $N_b \gg N_{\unknown}$, particularly for small diffusion. Thus the ratio of the cost of the particle method \eqref{eq:cost_particle} to the reference method \eqref{eq:cost_reference} is
\begin{equation}
  \frac{C_{particle}}{C_{reference}} = \bigO \left( \frac{C_{b}N_{u}}{N_{b}^2} \right), \quad\text{so typically } \frac{C_{particle}}{C_{reference}} \ll \bigO \left( \frac{C_{b}}{N_{b}} \right).
  \label{eq:cost_ratio}
\end{equation}
Thus for applications with a small enough number of observations or sufficient parallelism that a substantial proportion of the particles can be computed in parallel, the particle method should provide substantial speedup for large $N_{b}$ and in particular for problems with high-dimensional unknowns.

We have of course, made some simplifications in this analysis: For some choices of $\{\phi\}$ (e.g., discontinuous Galerkin \cite{hesthaven2007nodal}), the matrix $A$ may be sparse or block diagonal, reducing the cost of the matrix multiply and increasing the effectiveness of parallelization. However, we have also ignored some of the costs of \cref{alg:reference} and most of the performance improvement comes from not computing the full field $\pdesol$, so the computational complexity for large $N_\unknown$ or small $N_o$ should be as indicated here for a large class of reference algorithms.

The analysis for the boundary value method, \cref{alg:bvp_method}, is similar. The computational cost for the method is similar to that of \cref{alg:ad_method}; the main difference is that the number of timesteps required to reach the boundary of $\spatdom$ is uncertain and so may be larger or smaller than the largest $t_j$ for the time dependent problem. A Galerkin approach to solving \eqref{eq:bvp}, meanwhile, would replace the time integration of the reference method \cref{alg:reference} with an iterative solver of a linear system $A \pdesolvec = F$. Many efficient methods exist for solving these systems; however, because the number of basis functions required to represent $\pdesol$ would typically be much higher than the number of degrees of freedom in the parameter $\unknown$, we still expect the particle method to yield significant efficiencies for applications with sparse observations (i.e., small $N_o$).

\subsection{Limitations}
For completeness, we identify two key limitations of using \cref{alg:ad_method} or \cref{alg:bvp_method} to compute $\G$:
\begin{itemize}
  \item It is sometimes desirable to compute charateristics of $\pdesol$ beyond what is contained in $\G$. For example, \cite{borggaard2018bayesian} presents the statistics of the variance and variance dissipation of $\pdesol$ according to the Bayesian posterior on $\vfield$. Of course, because these particle methods do not involve computing the full field $\pdesol$, we cannot compute characteristics of $\pdesol$ beyond the values contained in $\G$. However, we note that the proposed particle algorithm could be used to compute the solution to the inverse problem, e.g., the maximum a posteriori $\unknown_{MAP}$. Then a single, computationally-expensive PDE solve could be used to compute the characteristics of that solution -- see, for example, the plot of $\pdesol(\unknown_{MAP})$ in \cref{fig:ex3:map} or the plot of the optimal $\pdesol$ in \cref{fig:bvp:ex1:results}.
  \item Many approaches to inverse problems (e.g., Hamiltonian Monte Carlo \cite{bou2018geometric,neal2011mcmc} for Bayesian inverse problems or Newton's method for optimization \cite{nocedal2006numerical}) require computing the \Frechet derivative $D\G(\unknown)$, which cannot be computed via the particle method in its current form. The gradient of $\G$ would therefore have to be approximated, e.g. via a finite-difference approximation involving multiple computations of $\G$. 
\end{itemize}

\section{Numerical Examples}\label{sec:results}
In this section, we provide examples demonstrating the applications and power of the particle methods \cref{alg:ad_method} and \cref{alg:bvp_method}. The examples were mostly written in the Julia numerical computing language \cite{bezanson2017julia}; the particle method was written in a combination of C and CUDA. The computations were run on the computational clusters at Virginia Tech,\footnote{\url{http://www.arc.vt.edu}} with each particle computation using a single NVIDIA P100 GPU. We note that a multi-GPU implementation should yield further speedup for larger problem sizes, with minimal overhead due to the embarrassingly parallel nature of the particle method.

\subsection{Example: Bayesian Inference of Fluid Flows}\label{sec:results:bayes}
In this section we present some numerical examples of the application of \cref{alg:ad_method} to the Bayesian inverse problem \cite{kaipio2005statistical,dashti2017bayesian} of estimating the background flow $\vfield$ from noisy point observations of $\pdesol$ with a constant diffusion coefficient:
\begin{equation}
  \frac{d\pdesol}{dt}(\x) = -\vfield(\x) \cdot \nabla \pdesol(t,\x) + \conductivity \Delta \pdesol(t,\x), \quad \nabla \cdot \vfield=0, \quad \pdesol(0,\x)=\pdesol_0(\x).
  \label{eq:adr}
\end{equation}
(In the notation of \eqref{eq:pde_diffusion}, $\sigma_{ij}(\x)=\sqrt{2\conductivity}\delta_{ij}$; note that since $\sigma$ is constant for this application, the Euler-Maruyama \eqref{eq:euler} and Milstein \eqref{eq:millstein} approximations are equivalent.) Accurate simulation of the low-$\conductivity$ case is an active area of research \cite{stynes2013numerical,morton1996numerical,codina2000stabilized,hundsdorfer2013numerical,mudunuru2016enforcing}; in this case, the solution $\pdesol$ can become high-dimensional very quickly, making accurate PDE solves highly computationally expensive \cite{borggaard2018bayesian}. Thus the particle method will yield substantial benefit as the dimension of the background flow increases. We represent the background flow $\vfield$ via a divergence-free Fourier expansion in terms of wave numbers $\kbf \in \Z^2$
\begin{equation}
	\vfield(\x) = \sum_{\kbf} \vfieldnobf_{\kbf} \frac{\kbf^\perp}{\norm{\kbf}} e^{2\pi i \kbf \cdot \x} 
  \label{eq:v_expand}
\end{equation}
where $\kbf^\perp = [-k_y, k_x]$ so that $\kbf \cdot \kbf^\perp = 0$, and $\vfieldnobf_{\kbf}$ obeys the reality condition $\vfieldnobf_{\kbf}=-\overline{\vfieldnobf_{-\kbf}}$. Thus, evaluating $\vfield(\x)$ involves computing a series of sines and cosines; this represents the dominant cost in the particle method.

We will now compute the Bayesian posterior distribution for three example problems while using \cref{alg:ad_method} to compute the forward map $\G$. All results are computed using the preconditioned Crank-Nicholson Markov Chain Monte Carlo (MCMC) method. For the details of this particular application, we refer the reader to \cite{borggaard2018bayesian}. 

The first example will demonstrate consistency of the method with the traditional two-step method implemented in \cite{borggaard2018bayesian}. The second will demonstrate how the particle method allows extension of Example 1 to higher-dimensional vector fields. Finally, we conclude with an example where \cref{alg:ad_method} allows application to a problem where the true background flow contains many thousands of components. We note that for these Bayesian inference experiments, for simplicity we have used $N_{b}=N_{\unknown}$ in the reference method; many applications would require $N_{b} \gg N_{\unknown}$ (more basis functions to represent $\pdesol$ than $\vfield$), increasing the cost of the reference method.

\subsubsection{Example 1: Consistency}
In \cite{borggaard2018bayesian}, the forward map $\G$ was computed with a two-step ``solve, then observe'' method: \eqref{eq:adr} was expanded in a Fourier basis, projected onto a system of ODEs, integrated in time, and then observed. In this section, we repeat the computations from Section 5.2 of that paper, in which $\vfield$ was assumed to have dimension of less than or equal to $197$ and $\G$ involved 100 point observations of $\pdesol$ (i.e., $N_{\unknown}=197$ and $N_o=100$). In this case, the dimension of $\vfield$ is low relative to the number of observations; thus the value of the particle method is limited -- the PDE solves in \cite{borggaard2018bayesian} took approximately the same time as the GPU particle implementation. The results are shown simply to demonstrate that they are consistent for the two methods. 

\begin{figure}[!htbp]
\centering
\begin{minipage}[b]{0.49\textwidth}
  \includegraphics[width=\textwidth]{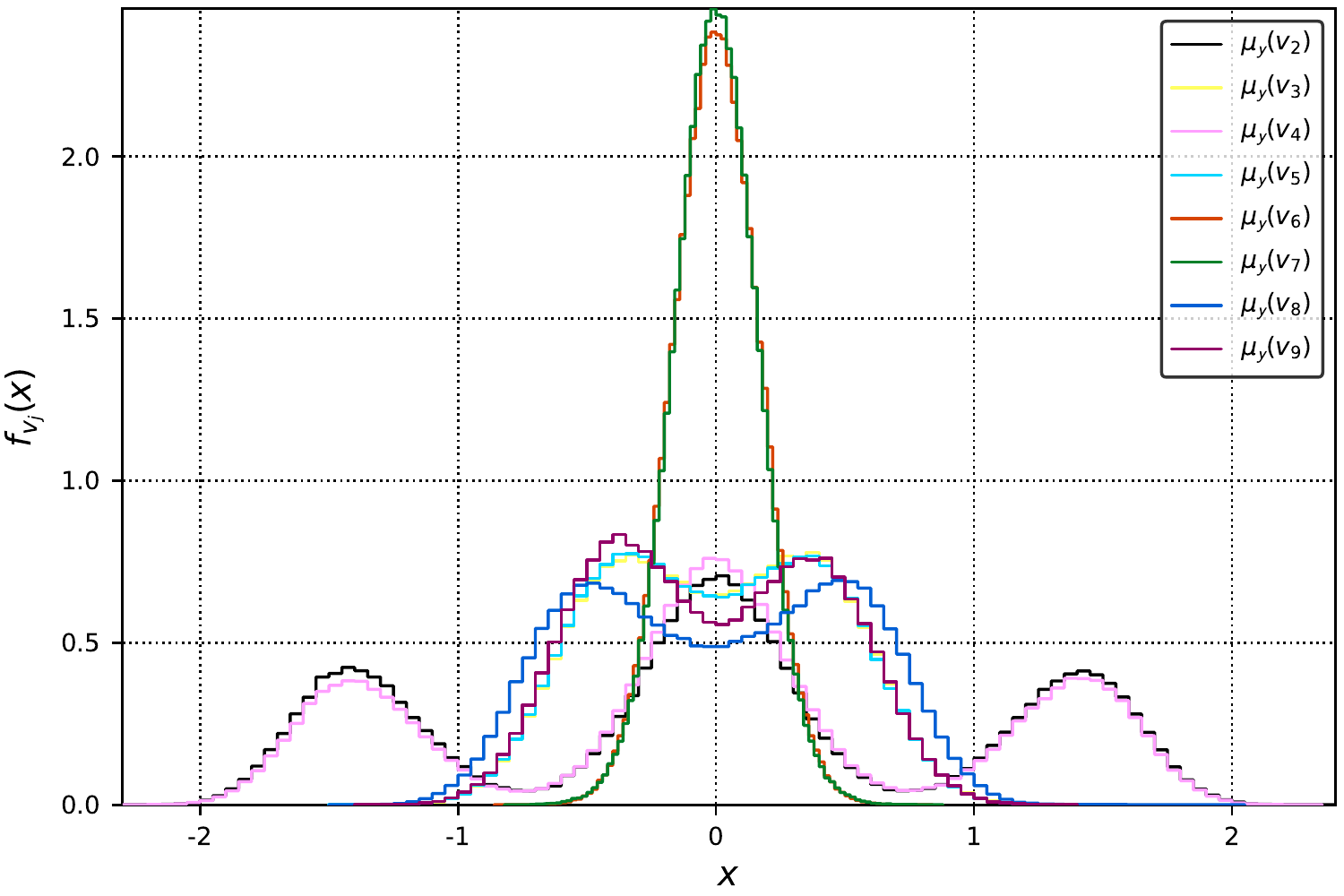}
\end{minipage}
\hfill
\begin{minipage}[b]{0.49\textwidth}
  \includegraphics[width=\textwidth]{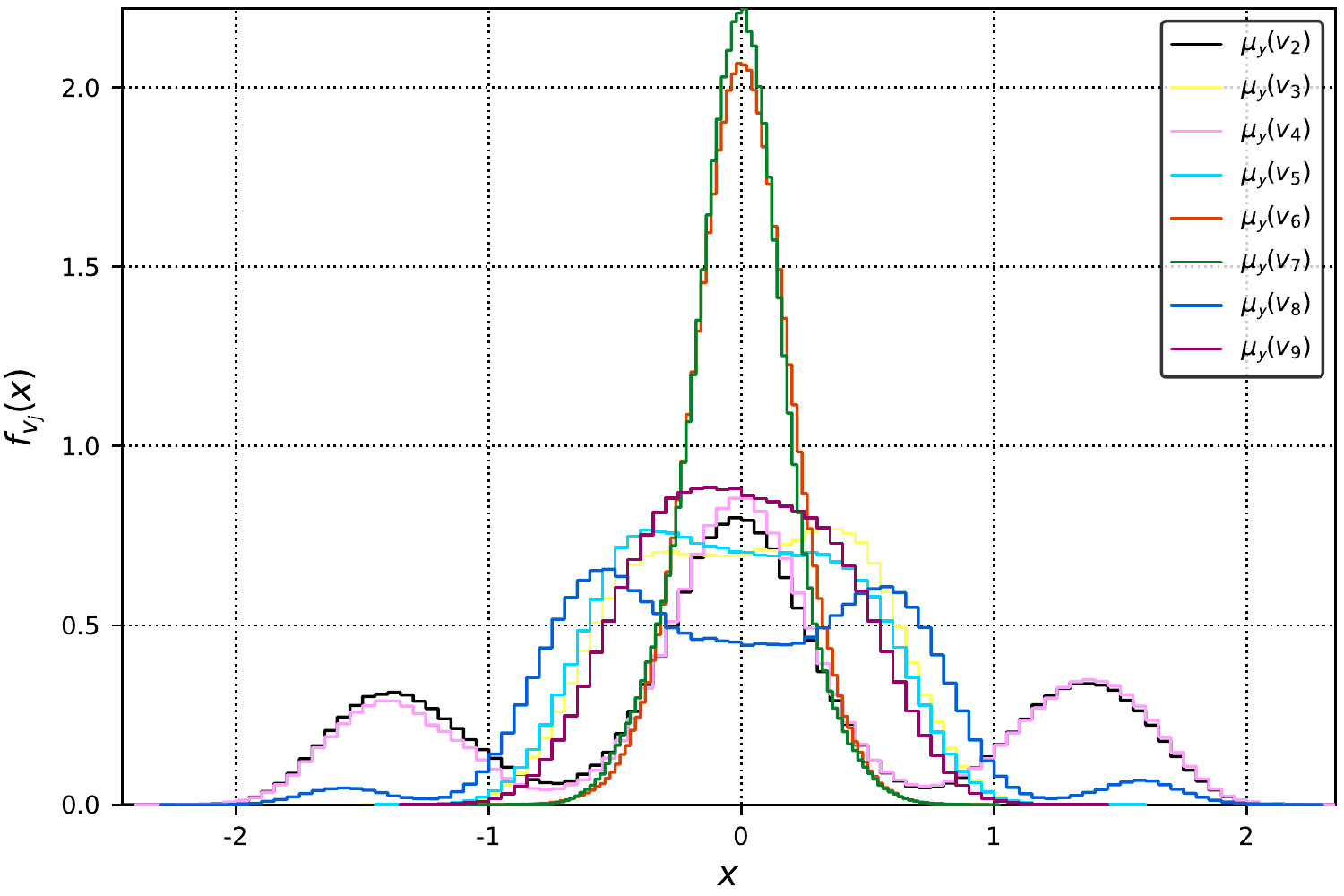}
\end{minipage}
\caption{\label{fig:ex1_prior_posterior} Posterior marginal distributions for the first eight components of $\vfield$. Left: 25M samples, Traditional Method, see \cite{borggaard2018bayesian}. Right: 10M samples, Particle Method.}
\end{figure}

%

\begin{figure}[!htbp]
\centering
\begin{minipage}[b]{0.49\textwidth}
  \includegraphics[width=\textwidth]{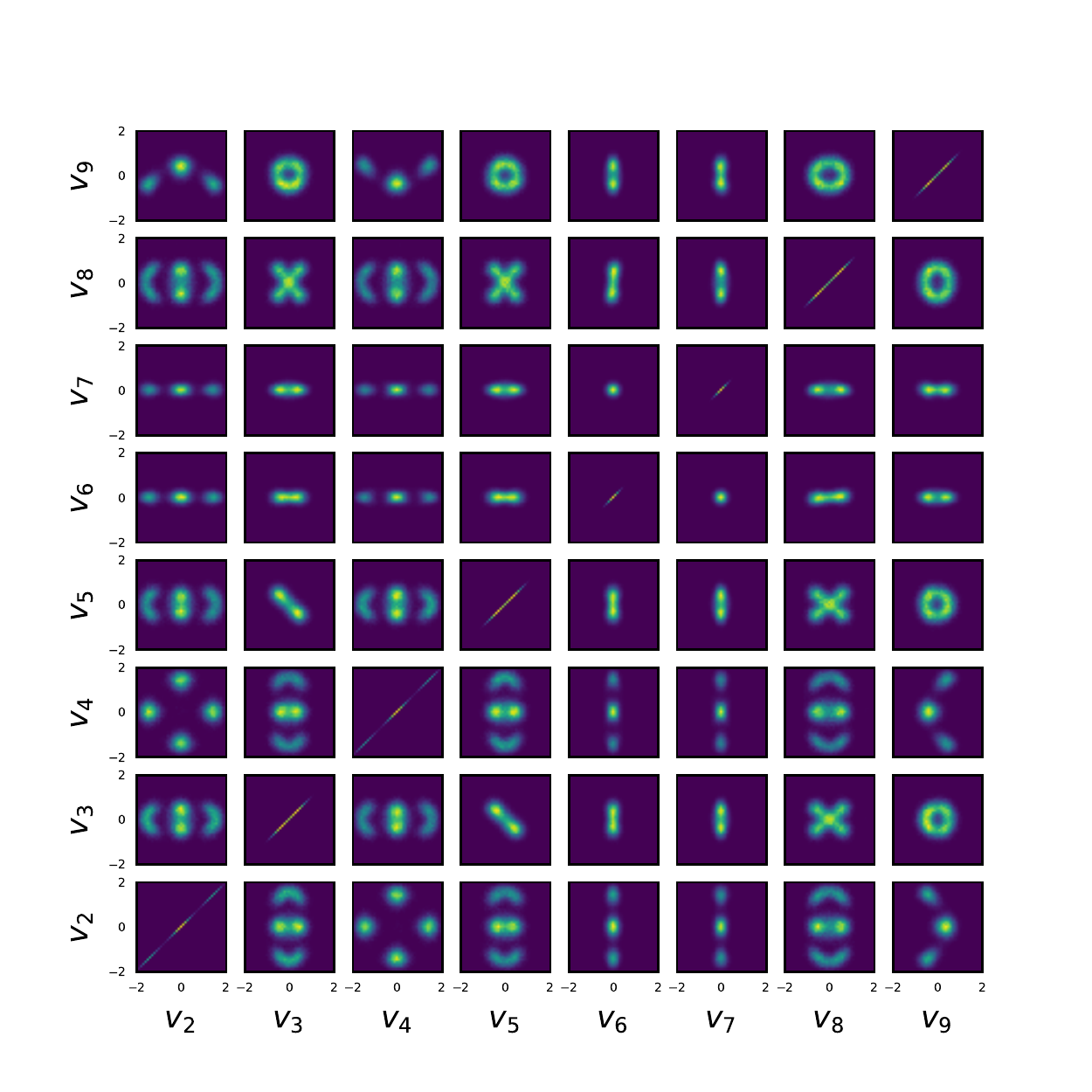}
\end{minipage}
\hfill
\begin{minipage}[b]{0.49\textwidth}
  \includegraphics[width=\textwidth]{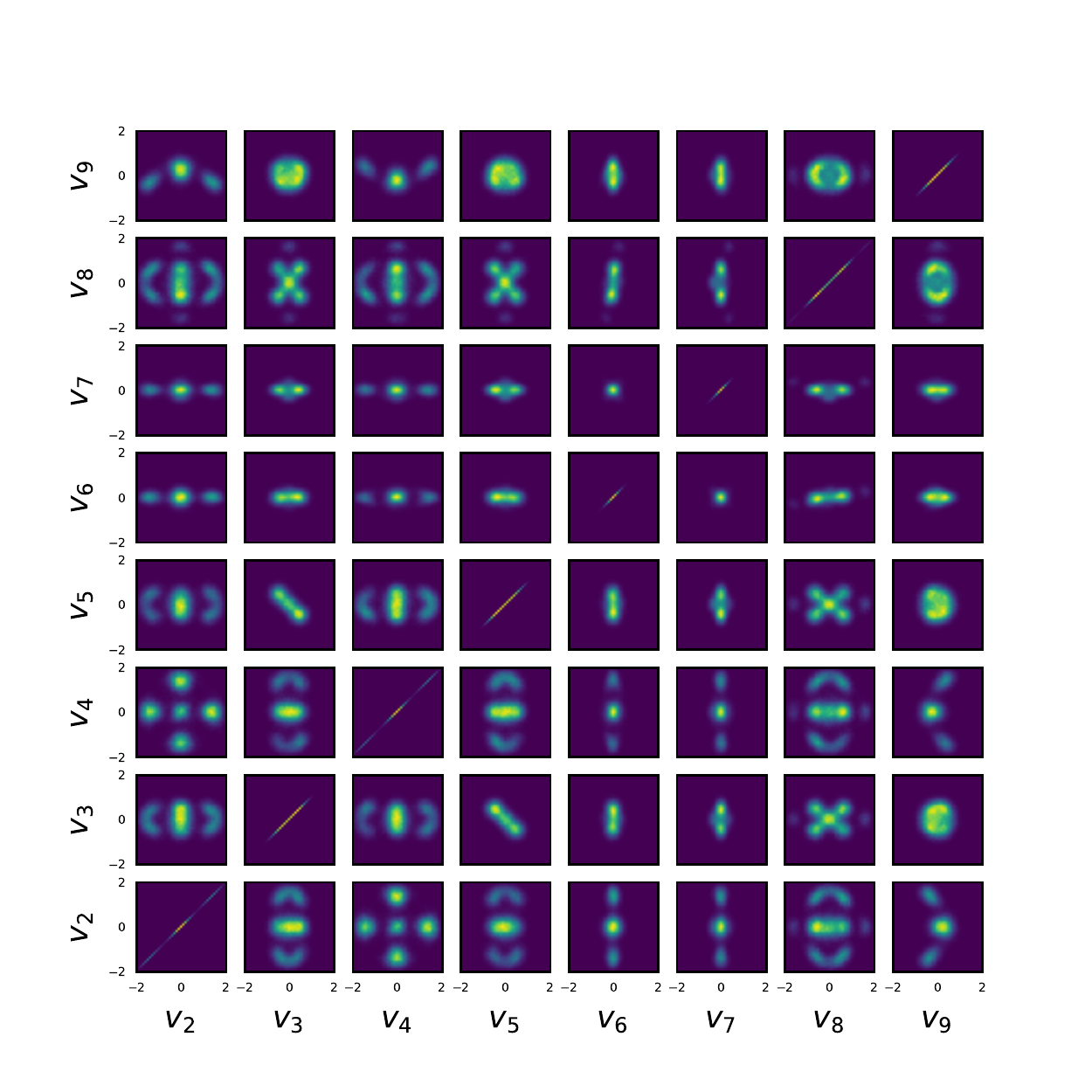}
\end{minipage}
\caption{\label{fig:ex1_hist2d}Posterior two-dimensional histograms for the first eight components of $\vfield$. Left: 25M samples, Traditional Method, see \cite{borggaard2018bayesian}. Right: 10M samples, Particle Method.}
\end{figure}

\cref{fig:ex1_prior_posterior} shows the computed posterior marginal histograms for the first eight Fourier components of the background flow $\vfield$, for both the original computation from \cite{borggaard2018bayesian} and as computed via \cref{alg:ad_method}. The structure of the posterior distributions have largely the same structure; smaller differences result from the fact that the posterior for this problem is highly complex and very slow to converge. \cref{fig:ex1_hist2d} shows the analogous two-dimensional histograms, which again show similar structure.

%

\subsubsection{Example 2: Extension to Higher-Dimensions}

In this example, we relax the assumption from Example 1 that the background flow has wave numbers with $\norm{\kbf}_2 \le 8$ ($N_{\unknown}=197$); here we allow $\norm{\kbf}_2 \le 32$ ($N_{\unknown}=3,209$). When computed using the reference method, the $16$-fold increase in the dimension of $\vfield$ yielded an increase of over $300$ in the computational cost of each sample, making it computationally intractable to generate enough samples to resolve the complex structure of the posterior distribution. By contrast, the computational cost of the GPU-based particle method only increased by a factor of $9$, as some of the linear computational cost \eqref{eq:cost_particle} was absorbed by the parallelism in the GPU. 

\begin{figure}[!htbp]
\centering
\begin{minipage}[c]{0.49\textwidth}
  \includegraphics[width=\textwidth]{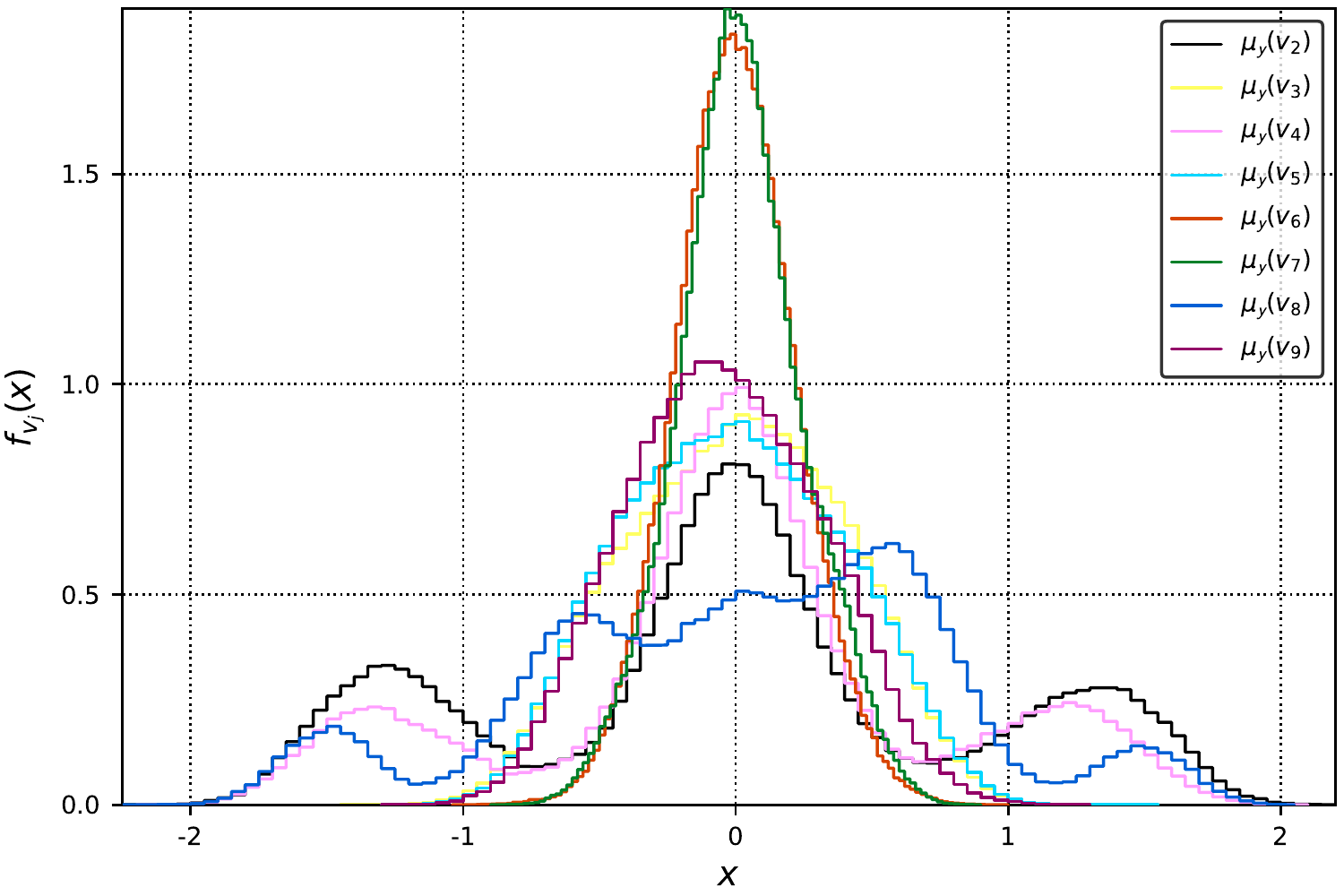}
\end{minipage}
\hfill
\begin{minipage}[c]{0.49\textwidth}
  \includegraphics[width=\textwidth]{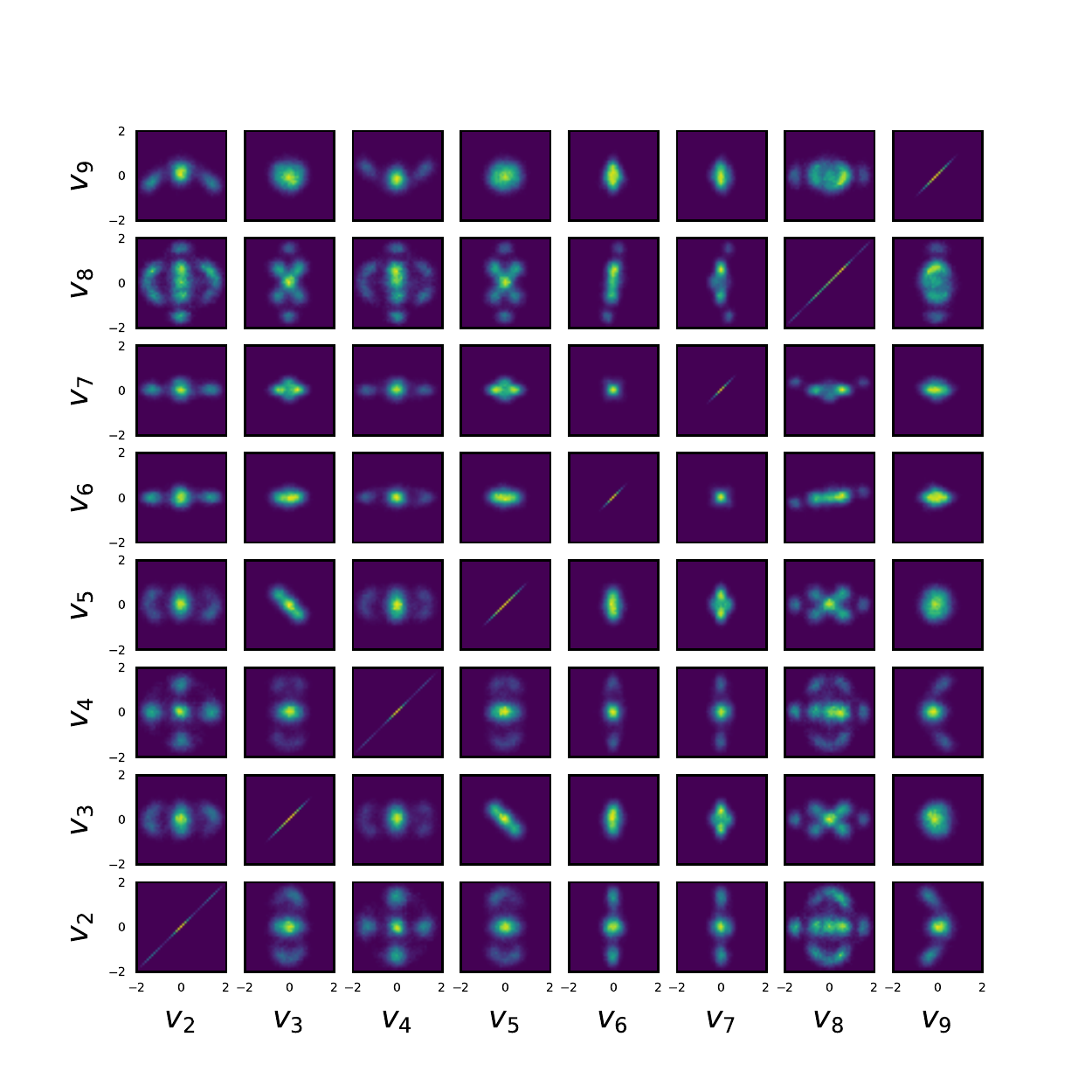}
\end{minipage}
\caption{\label{fig:ex2_posterior} Posterior marginal (left) and two-dimensional histograms for the first eight components of $\vfield$.}
\end{figure}

To approximate the posterior, we generated 2.5 million samples via 100 separate 25,000-sample chains. \cref{fig:ex2_posterior} shows the resulting posterior marginal distributions (left) and two-dimensional posterior histograms (right). When compared with \cref{fig:ex1_prior_posterior} and \cref{fig:ex1_hist2d}, respectively, we see that the posterior shows additional regions of probability mass - see, for example, the peaks for $\vfieldnobf_8 \approx \pm 1.5$. This indicates that the restriction of $\vfield$ to Fourier modes with $\norm{\kbf}_2 \le 8$ caused some possible candidate vector fields to be missed or de-emphasized.

%

\subsubsection{Example 3: Turbulent Background Flows} \label{sec:ex3:highd}
We conclude the Bayesian examples by using the particle method to address a problem of much higher dimension than is considered in \cite{borggaard2018bayesian}. We again consider the Bayesian inverse problem of estimating $\vfield$ from measurements of $\pdesol$ (see \eqref{eq:adr}). However, this time we consider $\vfield$ made up of components of wave numbers with $\norm{\kbf}_2 \le 80$, a total of 20,081 components ($N_{\unknown}=20,081$). We use the same initial condition $\pdesol_0$ and prior structure as in the previous two examples (though we allow the prior to extend to higher dimensions). In this example, we use $\conductivity=3 \times 10^{-5}$ from \cite{chen1998simulations}, which is typical for diffusion in water. We generate data by drawing a velocity field from the prior and computing $\pdesol$ at 13 evenly-spaced times between 0 and 0.5 at each of two observation locations: $[0,0], [\frac{1}{2},\frac{1}{2}]$. 
  
\begin{figure}[h]
\centering
\includegraphics[width=0.7\textwidth]{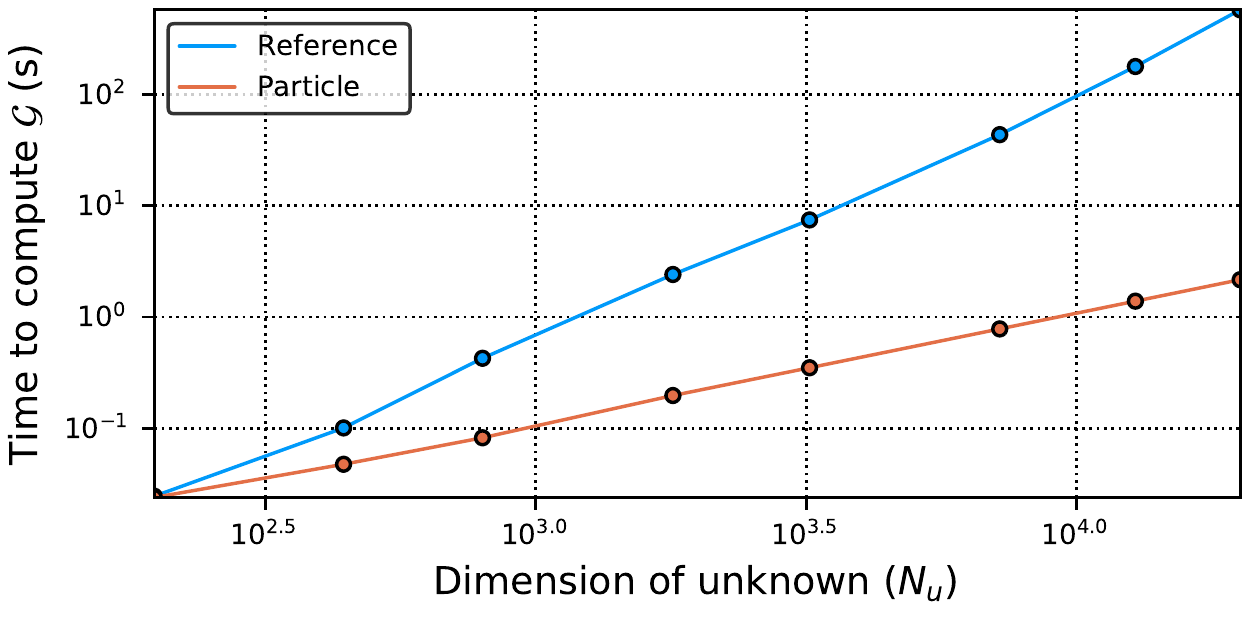}
\caption{\label{fig:gpu_time_compare} Time in seconds to compute $\G$ by dimension of unknown ($N_u$) for reference (Fourier) and particle methods.}
\end{figure}

Computing a Galerkin approximation of $\pdesol$ for this problem would require many matrix multiplications each of tens of thousands of rows and columns, making the computation of many thousands of samples, in general, computationally intractable. Using \cref{alg:ad_method} to compute the forward map $\G$, however, the computation scales roughly linearly with the number of unknowns, as shown in \cref{fig:gpu_time_compare}.

\begin{figure}[h]
\centering
\begin{minipage}[b]{0.49\textwidth}
  \includegraphics[width=\textwidth]{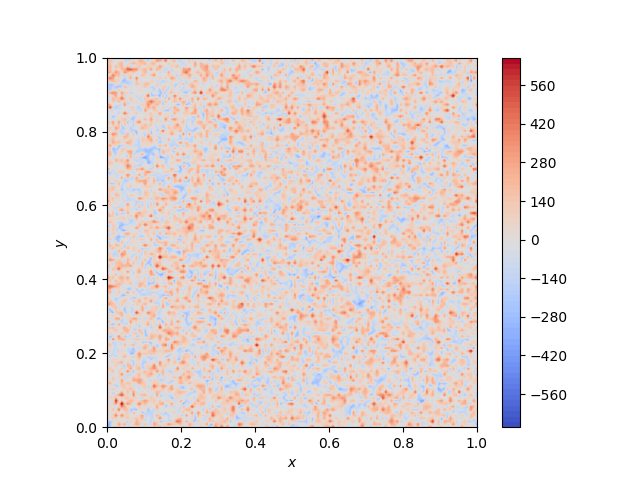}
\end{minipage}
\hfill
\begin{minipage}[b]{0.49\textwidth}
  \includegraphics[width=\textwidth]{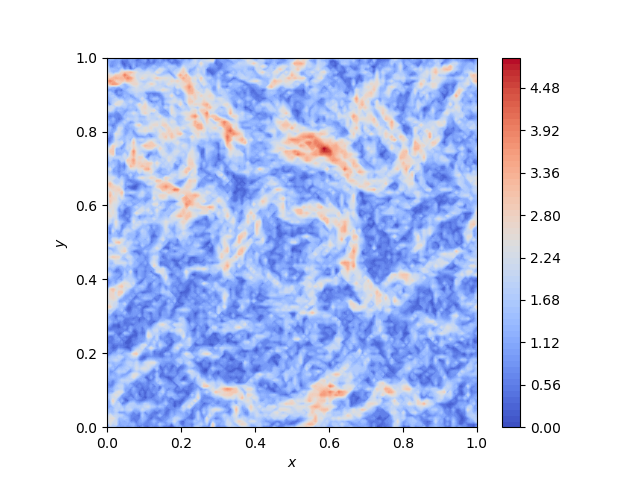}
\end{minipage}
\caption{\label{fig:ex3:map} Contour plots of vorticity (left) and norm (right) of $\vfield_{MAP}$.}
\end{figure}

\cref{fig:ex3:map} shows the vorticity and norm of the maximum a priori (MAP) point from 100,000 computed pCN MCMC samples; this can be thought of as the background flow that best matches both the prior measure and the observations for a given model of observational noise. The figures show the complexity of the background flows that could be considered in the inference by allowing the inclusion of tens of thousands of parameters.

\subsection{Optimal Forcing for Boundary Value Problem}
\begin{figure}[h]
  \centering
  \includegraphics[width=0.4\textwidth]{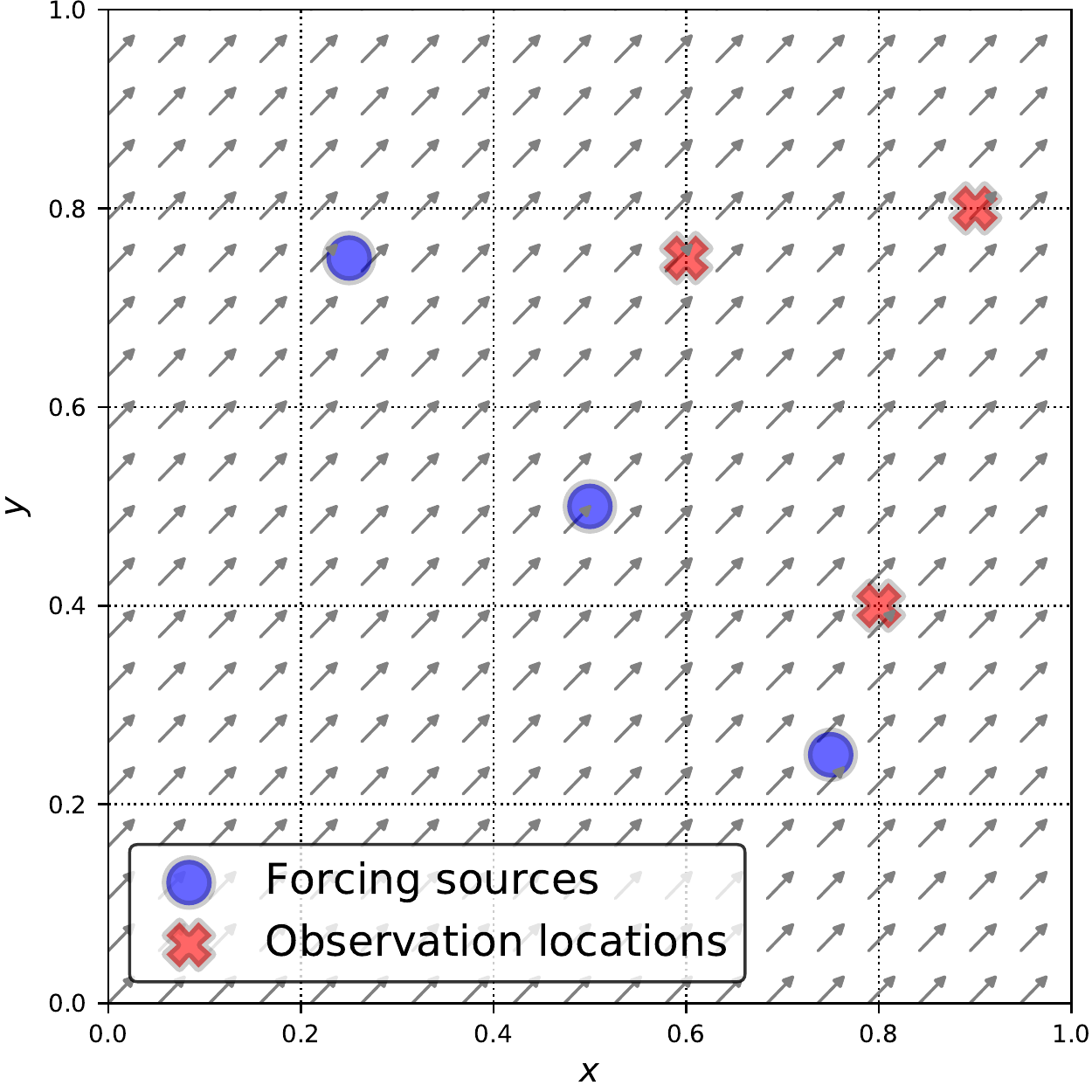}
  \caption{\label{fig:bvp:ex1:problem} Forcing locations, observation locations, and background flow.}
\end{figure}

In this section, we provide an example application of boundary value method \cref{alg:bvp_method}. For this example, we consider a version of the boundary value problem \cref{prob:bvp} with laminar background flow, constant diffusion, known boundary conditions, but to-be-determined forcing: 
\begin{equation}
  \begin{aligned}
    -\left[ \begin{array}{c} 1 \\ 1 \end{array}\right] \cdot \nabla \pdesol(\x) + \conductivity \Delta \pdesol(\x) &= \sum_{j=1}^3 f_j \exp \left( -4 \norm{\x-\x_{f}^{(j)}}_2^2 \right) && \quad \x \in \spatdom \coloneqq [0,1]^2 \\
    \pdesol(\x) &= \frac{1}{2} \left[ \cos\left( \frac{\pi}{2} x \right) + \cos\left( \frac{\pi}{2} y \right) \right] && \quad \x\in \partial\spatdom.
  \end{aligned}
  \label{eq:bvp:ex1}
\end{equation}
The goal of the problem is to find the forcing coefficients $\mathbf{F}=\left[ f_1, f_2, f_3 \right]$ that produce the scalar field $\pdesol$ that best matches the target values $\data$ at three observation locations $\left\{ \x_{obs}^{(i)} \right\}$. That is, we seek
\begin{equation}
  \mathbf{F}^\star \coloneqq \argmin_{\mathbf{F}} \norm{ \data - \G\left(\mathbf{F}\right) }, 
  \quad \G\left(\mathbf{F}\right) = \left[ \pdesol\left( \x_{obs}^{(1)},\mathbf{F} \right), \pdesol\left( \x_{obs}^{(2)},\mathbf{F} \right), \pdesol\left( \x_{obs}^{(3)},\mathbf{F} \right) \right].
  \label{eq:bvp:ex1:cost}
\end{equation}
The problem is thus an optimal control problem and can be interpreted as, e.g.,
\begin{itemize}
  \item Find the heat sources or sinks that produce the desired temperature at important locations in a room
  \item Find the forcing that minimizes the concentration of a chemical at key locations in a system (see, e.g., \cite{chang1992optimal,whiffen1993nonlinear} for examples with time-varying forcing)
\end{itemize}
The forcing locations $\left\{ \x_{f}^{(j)} \right\}$ and observation locations $\left\{ \x_{obs}^{(i)} \right\}$ used in this example are shown in \cref{fig:bvp:ex1:problem}. For this problem, we use diffusion coefficient $\conductivity=0.282$, for water diffusing in air \cite{cussler2009diffusion}, and seek to match data $\data=[0, 0, 0]$.

Recall that the forward map $\G$ is sparse for this problem - we only need $\pdesol$ evaluated at three points. We can therefore leverage \cref{alg:bvp_method} to speed up the computations. We use a Nelder-Mead simplex method \cite{nelder1965simplex,lagarias1998convergence} as implemented in Julia's Optim package \cite{mogensen2018optim} to seek the optimal $\mathbf{F}$. The method required evaluating $\G$ for thousands of possible values of $\mathbf{F}$. The requirement of multiple evaluations of the forward map is typical of many approaches to PDE-constrained optimization problems \cite{hinze2009optimization}. The choice of Nelder-Mead is merely used to demonstrate the effectiveness of \cref{alg:bvp_method} within an optimization setting.

\begin{figure}[h]
\centering
\begin{minipage}[b]{0.49\textwidth}
  \includegraphics[width=\textwidth]{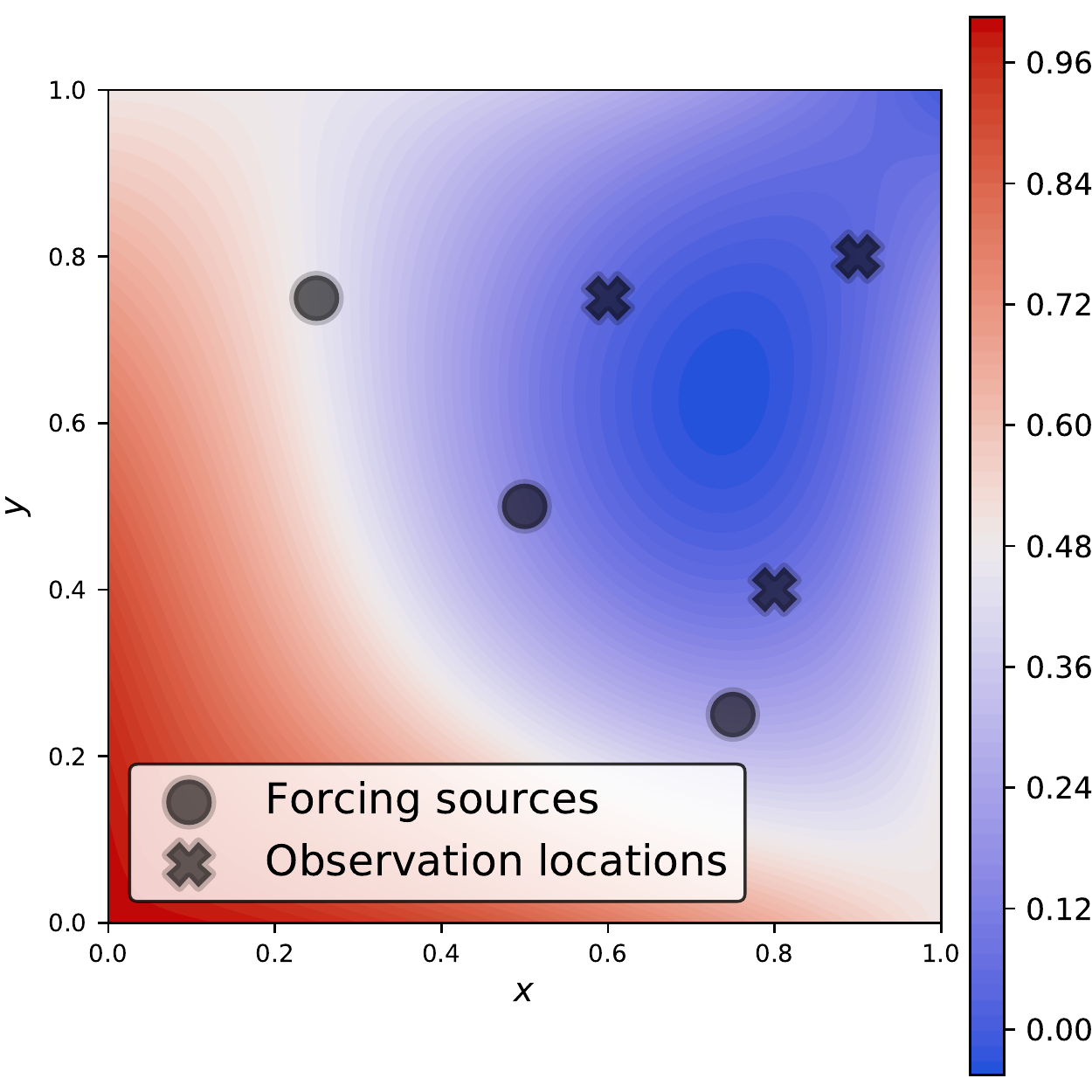}
\end{minipage}
\hfill
\begin{minipage}[b]{0.49\textwidth}
  \includegraphics[width=\textwidth]{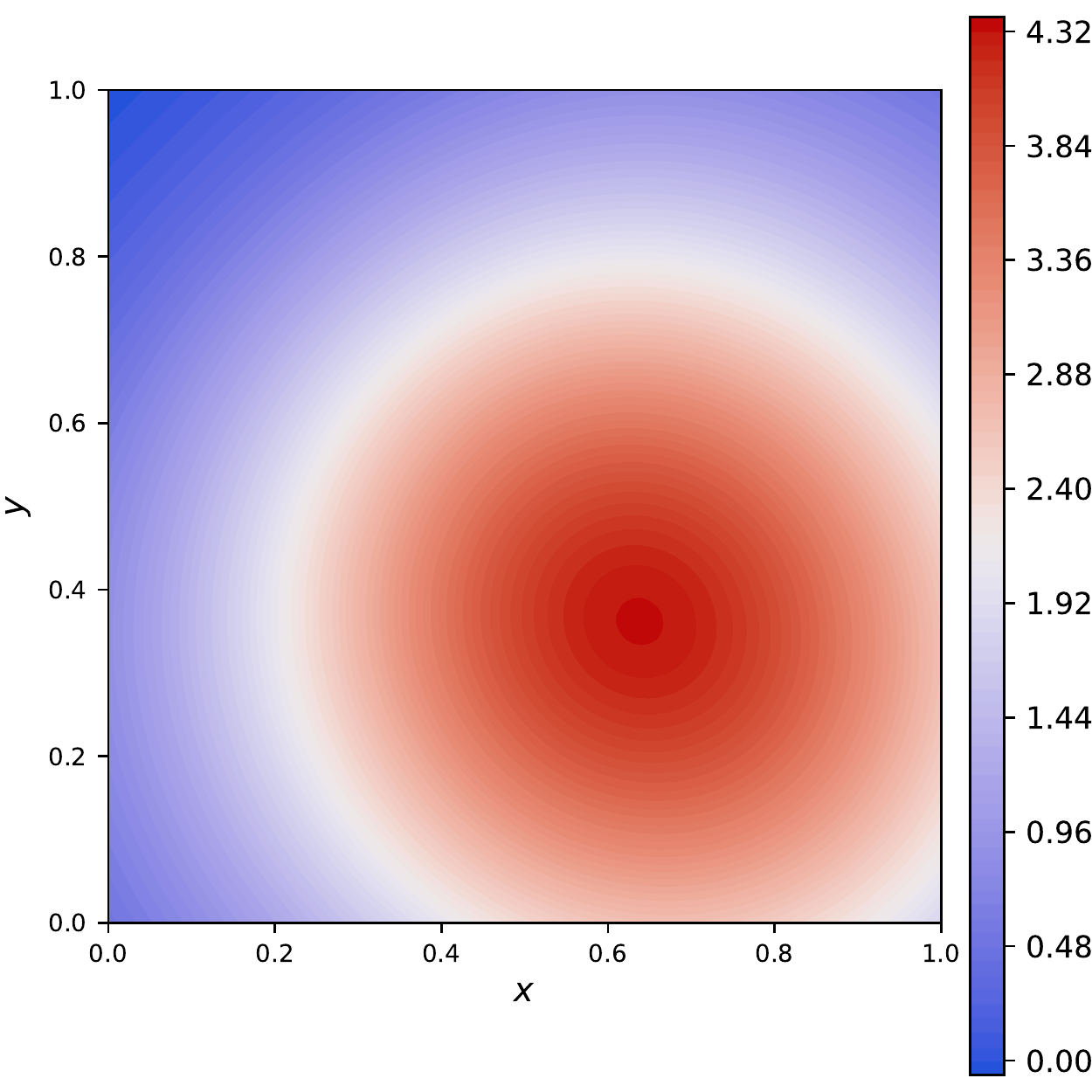}
\end{minipage}
\caption{\label{fig:bvp:ex1:results} Scalar field $\pdesol$ (left) and forcing function (right) for computed optimal coefficients $\mathbf{F}^\star$.}
\end{figure}

The results, $\pdesol$ and $f$, for the computed optimal forcing are shown in \cref{fig:bvp:ex1:results}; notice that the observation points lie along the contour $\pdesol=0$, indicating that we have found a set of parameters that are a good match for the data. (These plots, which required approximation of $\pdesol$ for the full domain, were generated via a finite element solver after $\mathbf{F}^\star$ was computed via the particle method; the plots required more work to generate than finding the optimal set of coefficients.)

We note that other, similar optimization problems -- for example, fixing the forcing and instead seeking $\vfield$, $\conductivity$, and/or boundary conditions that produce $\pdesol$ best matching the data -- could be addressed via an analogous approach, again by leveraging \cref{alg:bvp_method} to compute $\pdesol$ at the observation points for each value of the parameter considered by the optimization algorithm.

\section{Acknowledgments}
This work was supported in part by the National Science Foundation
under grants DMS-1313272, DMS-1522616, and DMS-1819110; the National Institute for
Occupational Safety and Health under grant 200-2014-59669; and the
Simons Foundation under grant 515990. We would also like to thank the
Mathematical Sciences Research Institute, the Tulane University Math
Department, and the International Centre for Mathematical Sciences
(ICMS) where significant portions of this project were developed and
carried out.

The authors acknowledge Advanced Research Computing at Virginia Tech\footnote{\url{http://www.arc.vt.edu}} for providing computational resources and technical support that have contributed to the results reported within this paper.

\addcontentsline{toc}{section}{References}
\begin{footnotesize}
\bibliographystyle{plain}
\bibliography{references}
\end{footnotesize}

\newpage

\vspace{.3in}
\begin{multicols}{2}
\noindent
Jeff Borggaard\\
{\footnotesize Department of Mathematics\\
Virginia Tech\\
Web: \url{https://www.math.vt.edu/people/jborggaa/}\\
Email: \href{mailto:jborggaard@vt.edu}{\nolinkurl{jborggaard@vt.edu}}} \\[.25cm]
Nathan Glatt-Holtz\\ {\footnotesize
Department of Mathematics\\
Tulane University\\
Web: \url{http://www.math.tulane.edu/~negh/}\\
Email: \href{mailto:negh@tulane.edu}{\nolinkurl{negh@tulane.edu}}} \\[.2cm]

\columnbreak

 \noindent Justin Krometis\\
{\footnotesize
Advanced Research Computing\\
Virginia Tech\\
Web: \url{https://www.arc.vt.edu/justin-krometis/}\\
Email: \href{mailto:jkrometis@vt.edu}{\nolinkurl{jkrometis@vt.edu}}} \\[.2cm]
 \end{multicols}

\end{document}